\documentclass[12pt]{article}
\usepackage{epsfig,amssymb,epic,eepic,color,graphics,amsmath}
\usepackage[utf8]{inputenc}
\usepackage[english]{babel}
\input{amssym.def}

\title{Existence of an energy function for 3-dimensional chaotic ``sink-source'' cascades}
\author{M. Barinova, V. Grines, O. Pochinka, B. Yu}
\date{DSA Laboratory HSE (Russia), Tongji University (China)}

\usepackage{graphicx}
\newtheorem{theo}{Theorem}

\newtheorem{defi}{Definition}
\newtheorem{lemm}{Lemma}

\newenvironment{demo}{{\bf Proof: }}{\hfill $\diamond$\medskip}

\begin{document}
\sloppy
\maketitle
\begin{abstract}
The paper is a continuation of research in the direction of energy function (a smooth Lyapunov function whose set of critical points coincides with the chain recurrent set of a system) construction for discrete dynamical systems. The authors established the existence of an energy function for any $A$-diffeomorphism of a three-dimensional closed orientable manifold whose non-wandering set consists of chaotic one-dimensional attractor and repeller.
\end{abstract}

\maketitle 

\section{Introduction and  statement of results}

Let $ f $ be a diffeomorphism of a smooth closed $ n $-manifold $ M^n $. One says that $f$ {\it satisfies  axiom A} if its non-wandering set is hyperbolic and the periodic points are dense in it. We say in this case that $f$ is an $A$-diffeomorphism. For $ A $-diffeomorphisms, the Smale spectral decomposition theorem  \cite{Smale1961} holds, that is non-wandering set is a  union of finite number of pairwise disjoint sets called {\it basic sets}, each of which is compact, invariant and topologically transitive.

A basic set $\Lambda $ is called an {\it attractor} of an A-diffeomorphism $ f $ if it has a compact neighborhood $ U_\Lambda $ such that $ f (U_\Lambda) \subset \textit {int} ~ U_\Lambda $ and 
$ \Lambda = \bigcap \limits_{k \geqslant 0} f^k {(U_\Lambda)} $.  $ U_\Lambda $ is called a {\it trapping neighborhood} of $\Lambda$. A {\it repeller} is defined as the  attractor for $ f^{- 1} $. By a {\it dimension} of the attractor (repeller) we mean its topological dimension. The set $\bigcup \limits_{k \in\mathbb Z} f^k {(U_\Lambda)}$ is called a {\it basin} of the attractor $ \Lambda $.

By Theorem 3 in \cite{Plykin1971}, every one-dimensional basic set  of an A-diffeomorphism on a surface is an attractor or repeller, and $ \dim {W_x^u} = \dim {W^s_x} = 1 $, for $ x \in \Lambda $. The trapping neighborhood for such basic set is a surface with a boundary (see section \ref{bss}). Therefore it can be  naturally included to a three-dimensional dynamics.

\begin{defi}\label{sa} A connected one-dimensional attractor of $ A $-diffeomorphism $ f: M^3 \to M^3 $ is called a {\rm canonically embedded surface attractor} if:
\begin{itemize}
\item $A$ has a trapping  neighborhood $ U_{A} $ of the form $ S_A \times [-1,1] $, where $ S_A = S_A \times \{0 \} $ is a surface with a boundary containing $ A $;
\item $S_A$ is a trapping neighborhood of $A$ as an attractor for the diffeomorphism $ g = f | _{S_A}$;
\item diffeomorphism $ f | _{U_A} $ is topologically conjugate to the diffeomorphism $ \phi (w, z) = (g (w), z/2),\, (w, z) \in S_A \times [-1,1 ] $.
\end{itemize}
\end {defi}
A one-dimensional repeller is called a {\it canonically embedded surface repeller} if it is a canonically embedded connected one-dimensional surface attractor for the diffeomorphism $ f^{-1}$.

One naturally expect to construct a diffeomorphism on a closed 3-manifold by a gluing dynamics from canonically embedded surface one-dimensional attractor and repeller. In this paper the problem was solved (see section \ref{gl}) for the attractor and repeller models obtained by DA-surgeries from the same Anosov diffeomorphism. It is clear that the construction admits generalizations to the attractor and repeller obtained by DPA-surgeries from the same pseudo-Anosov diffeomorphism. 

\begin{theo}\label{theo0} There are infinitely many pairwise $\Omega$-non-conjugated $3$-diffeomorphisms whose non-wandering set is a union of canonically embedded one-dimensional surface attractor and repeller.
\end{theo}

Note that the similar ``one-dimensional attractor-repeller'' dynamics on the surface also always is not structurally stable due to  results by R. Robinson and R. Williams \cite{RoWi}. The  construction of   3-diffeomorphisms with one-dimensional attractor-repeller (not surface) dynamics firstly was suggested in \cite{Gib72}, all examples also were not structurally stable. 
In \cite{Bonatti2010}, \cite{Shi2014} structurally stable examples with one-dimensional attractor-repeller  dynamics were constructed. But the question  whether the basic sets are  canonically embedded in  surfaces are open because the construction is very different from one suggested in this paper.
The surface dynamics allows us to prove the existence of an energy function for the examples constructed in this paper. 

{\it Lyapunov function} is  a continuous function that is constant on the chain components  and decreases along trajectories outside them. It is an important characteristic of a dynamical system, initially introduced by A.M. Lyapunov to study the stability of equilibriums of  differential equations systems. In modern dynamics it plays a significant role. The fundamental theorem of dynamical systems, proved by K. Conley \cite{Conley1978} in 1978, establishes the existence of a Lyapunov function for every dynamical system. 

If the Lyapunov function is smooth and the set of its critical points coincides with the chain recurrent set of the dynamical system, then it is called an {\it energy function}. In this case, many characteristics of a dynamical system directly follow from the properties of its energy function. Unlike flows, diffeomorphisms do not always have an energy function. Moreover, counterexamples are known even for cascades with a regular dynamics. 

First such example was constructed by D. Pixton \cite{Pixton1977} in 1977, it was a Morse-Smale diffeomorphism on 3-sphere. Necessary and sufficient conditions for the existence of an energy function for arbitrary Morse-Smale 3-diffeomorphism were found by V.Z. Grines, F. Laudenbach, and O.V. Pochinka in \cite{GrLauPo2012}.

For systems with chaotic dynamics first constructions of energy functions were done in  \cite{GrNoPo2015-1}, \cite{GrNoPo2015-2} for $ A $-diffeomorphisms with basic sets of co-dimension one on 2- and 3-manifolds. This paper is a  continuation of those papers for 3-diffeomorphisms with one-dimensional basic sets. 

\begin{theo} \label{theo1} Every $ A $-diffeomorphism of a closed orientable 3-manifold $M^3$, whose non-wandering set is a union of connected canonically embedded one-dimensional surface attractor and repeller, has an energy function.
\end{theo}

\section{One dimensional basic sets for diffeomorphisms of surfaces}\label{bss}
In this section we give a brief description of one dimensional basic sets for diffeomorphisms of surfaces. For simplicity everywhere below we assume that $ \Lambda $ is a connected attractor. 

By theorem 1 in \cite{Plykin1971},  $ \Lambda = \bigcup \limits_{x \in \Lambda} W^u_x $. In addition, by lemmas 2.1, 2.4, 2.5 in  \cite{Grines1975} at least one of the connected components of the set $ W^s_x \setminus \{x \}, x \in \Lambda $ contains a dense set in $ \Lambda $ and there are a finite number of points $ x \in \Lambda $ for which one of the connected components $W^{s\varnothing}_x$ of the set $ W^s_x \setminus \{x \} $ does not intersect $ \Lambda $. Such points are called {\it $ s $-boundary}, their set is not empty and consists of periodic points. The set $W^s_{\Lambda}\setminus \Lambda$ consists of a finite number path-connected components. 

\begin{figure}[htp]
\centerline{\includegraphics[width=8 true cm]{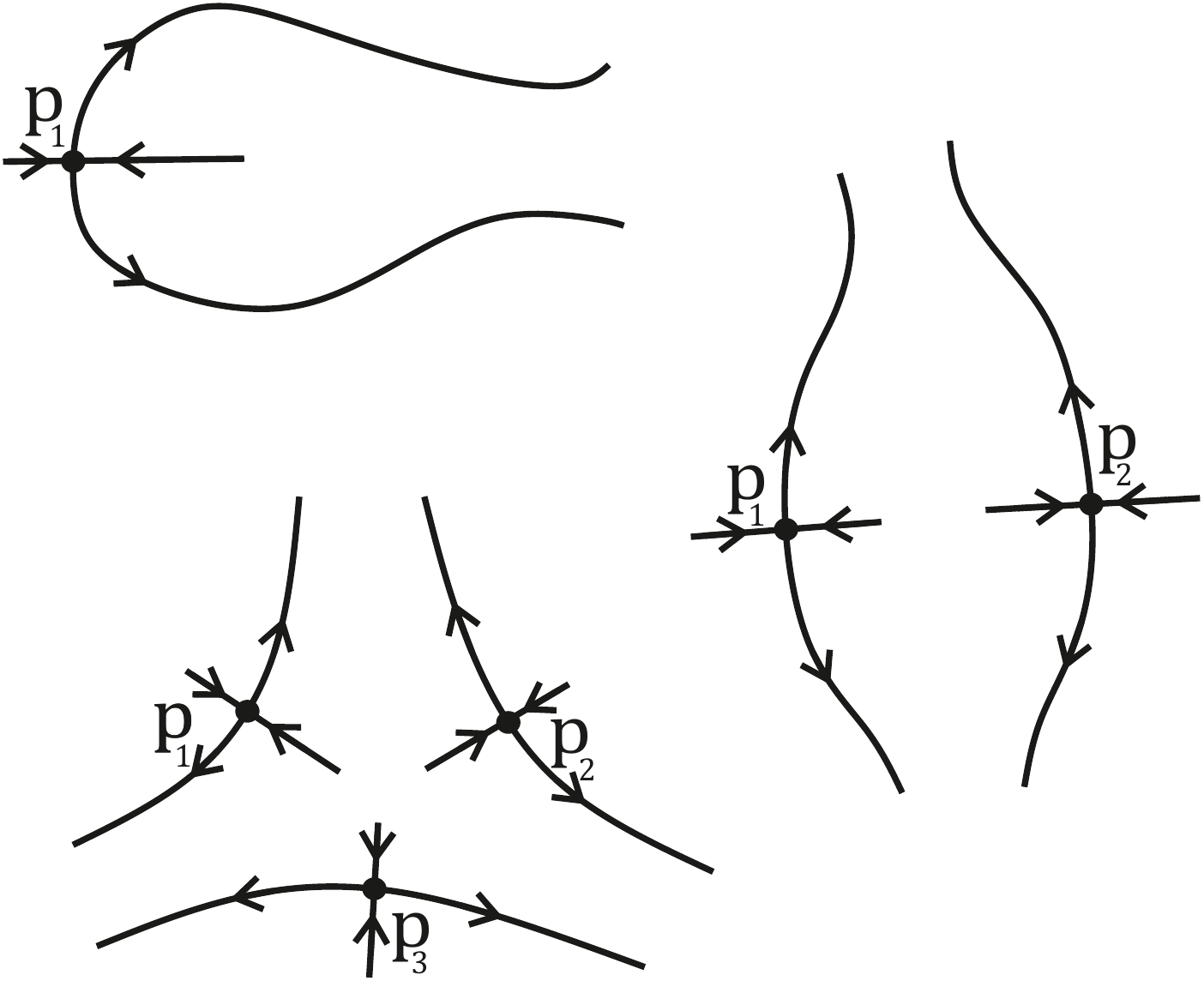}}
\caption{Bunches of degrees 1, 2 and 3}\label{periodical_points}
\end{figure}

{\it Bunch} $b$ of the attractor $\Lambda$ is the union of the unstable manifolds $W^u_{p_1}, \ldots, ~ W^u_{p_{r_b}}$ of $s$-boundary points $p_1, \ldots, ~ p_{r_b}$ of $\Lambda$ for which $W^{s\varnothing}_{p_1},\dots,W^{s\varnothing}_{p_{r_b}}$ belong to the same path-connected components of $W^s_{\Lambda}\setminus \Lambda$. The number $ r_b $ is called a {\it degree of the bunch} (see figure \ref{periodical_points}). Let $ B _{\Lambda} $ be the set of all bunches of the attractor $ {\Lambda} $. 

Let us explain how to construct a  trapping neighborhood for $\Lambda$. A proof of the existence of a neighborhood with properties below is possible to find, for example, in \cite{GrPo}. 

The property $\dim\,W^s_x=\dim\,W^u_x=1,\,x\in\Lambda$ allows  to introduce the notation $ (y,z)^s \, ((y,z)^u) $ for the arc of a stable (unstable) manifold  bounded by points $y,z$. Then for the bunch $ b \in B _{\Lambda} $ we choose a sequence of points $ x_1, \dots, x_{2r_b} $ such that:

\begin{itemize}

\item $x_{2j-1}$, $x_{2j}$ belong to different connected components of the set $W_{p_j}^u \setminus p_j$;

\item $x_{2j+1}\in W^s_{x_{2j}}$ (suppose $x_{2r_b +1} = x_1$);

\item $(x_{2j}, x_{2j+1})^s \cap\Lambda = \emptyset$, $j = 1, \dots, r_b$.

\end{itemize}

For $j \in \{1, \dots, r_b \}$ we choose a pair of points $ \tilde {x}_{2j-1} $, $ \tilde {x}_{2j} $, and a simple curve $ l_j $ with boundary points $ \tilde {x}_{2j-1} $, $ \tilde {x}_{2j} $ such that:

\begin{itemize}

\item $ (\tilde {x} _{2j}, \tilde {x} _{2j + 1})^s \subset (x_{2j}, x_{2j + 1})^s $;

\item the curve $ l_j $ intersects transversally with a stable manifold of any point of the arc $ (x_{2j-1}, x_{2j})^u $ exactly at one point;

\item the curve $ L_b = \bigcup \limits_{j = 1}^{r_b} [l_j \cup (\tilde {x} _{2j}, \tilde {x} _{2j + 1})^s] $ is a simple closed smooth curve and the set $ L _{\lambda} = \bigcup \limits_{b \in B _{\lambda}} L_b $ has the following properties:

\begin{itemize}

\item $ f (L _{\lambda}) \cap L _{\lambda} = \emptyset $;

\item for every curve $ L_b $, $ b \in B _{\lambda} $ there exists a curve from the set $ f (L _{\lambda}) $ such that these two curves are the boundary of the two-dimensional annulus $ K_b $;

\item annulus $ \{K_b, ~ b \in B_{\Lambda} \} $ are pairwise disjoint (see figure \ref{neighborhood_construction}).

\begin{figure}[htp]
\centerline{\includegraphics[width=8 true cm]{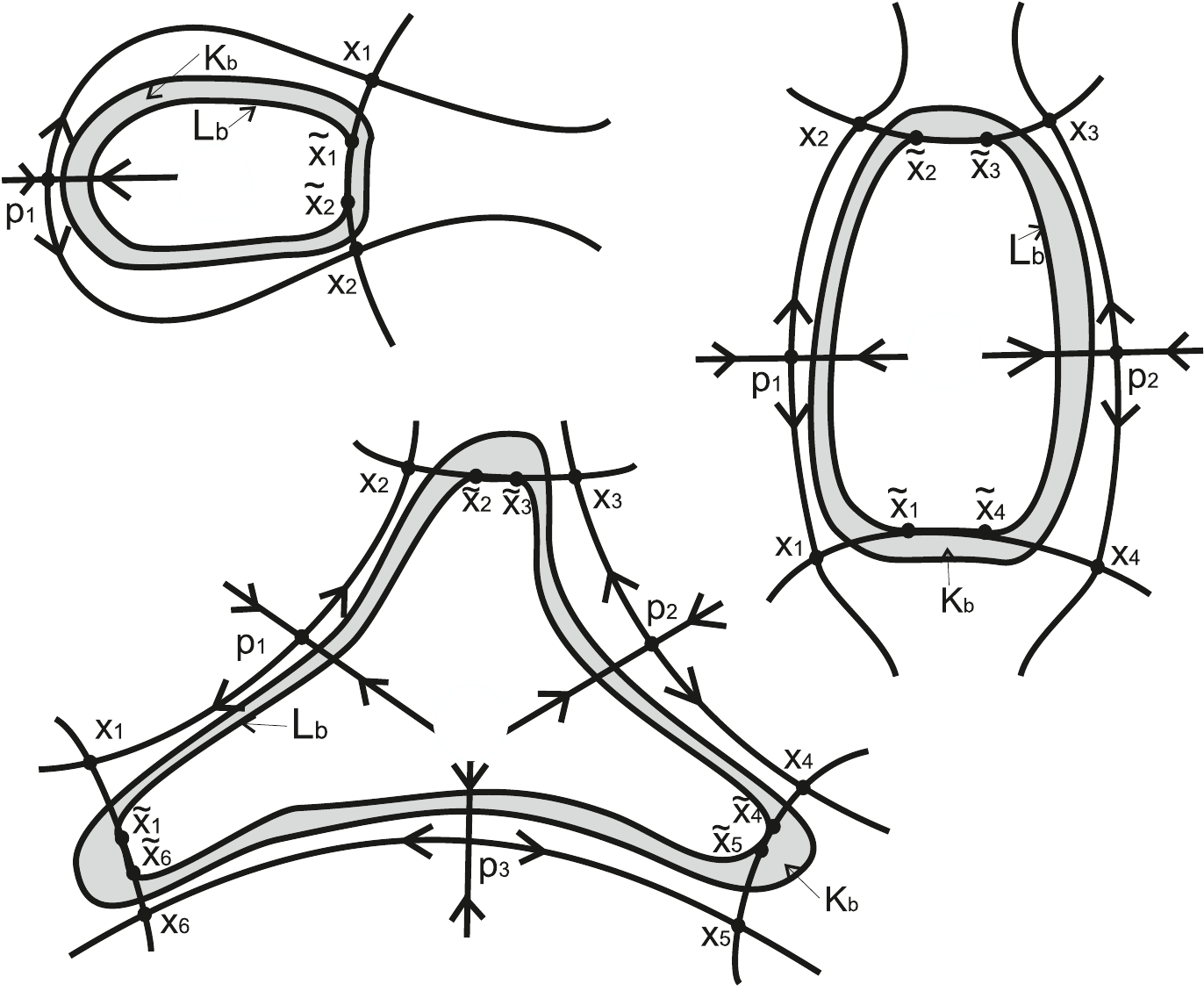}}
\caption{Construction of a  trapping neighborhood for $\Lambda$}\label{neighborhood_construction}
\end{figure}

\end{itemize}
\end{itemize}

Let $ K_\Lambda = \bigcup \limits_{b \in B_{\Lambda}} K_b $. It is directly verified that a surface with boundary $ U_\Lambda = \bigcup \limits_{k \in \mathbb{N}} f^k (K_\Lambda) \cup \Lambda $ is a trapping neighborhood of the attractor $ \Lambda $.
\section{Construction of canonically embedded surface  attractor (repeller)}\label{example}

\subsection{Anosov diffeomorphism on a 2-torus}

Let $C \in GL (2, \mathbb Z)$ be a hyperbolic matrix with eigenvalues $\lambda_1, \lambda_2$ such that $ \lambda = | \lambda_1 |> 1 $ and $ | \lambda_2 | = 1 / \lambda $. Since the matrix $ C $ has a determinant equal to 1, it induces a hyperbolic automorphism $ \widehat C: \mathbb T^2 \to \mathbb T^2 $ with a fixed point $ O $. This diffeomorphism is an Anosov diffeomorphism possessing two transversal invariant
foliations (stable and unstable) whose leaves are irratioanal windings on the torus.
In addition, the set of periodic points of the diffeomorphism $ \widehat C $ is also dense
on $ \mathbb T^2 $.

\subsection{Smale surgery} \label{surgery}

\begin{figure} [h]
\centerline{\includegraphics [width = 6 true cm]{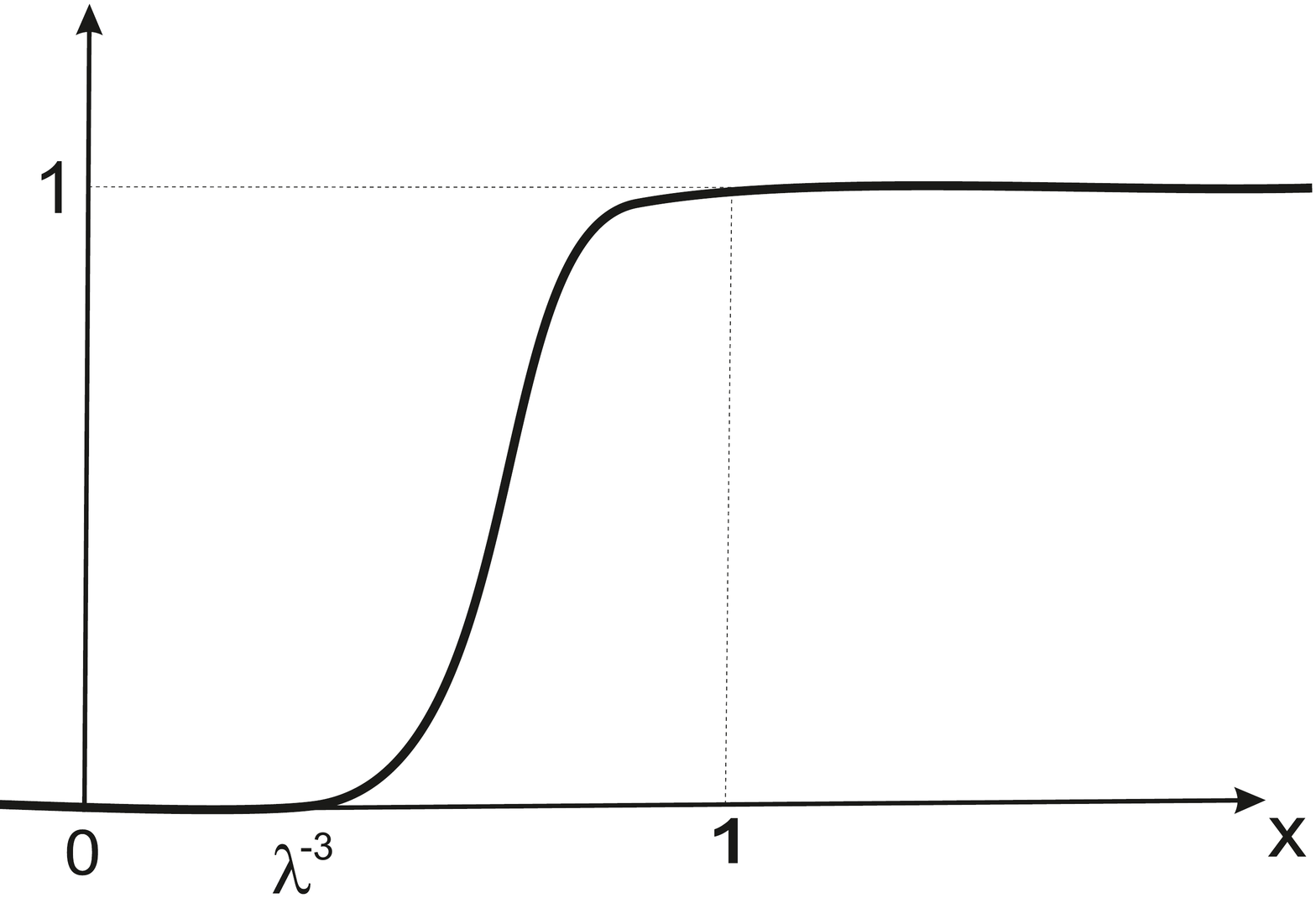}}
\caption{Graph of the function $ \sigma (x) $} \label {si}
\end{figure}

Let $ \sigma: \mathbb R \to [0,1] $ be a sigmoid defined by the formula (see the figure \ref{si})
$$ \sigma (x) = \begin{cases} 0, & x \leqslant \lambda^{- 3}, \\ \frac {1} {1+ \exp \left(\frac{\frac{\lambda^{- 3} +1} {2} -x} {{\left(x- \lambda^{- 3} \right)}^2 {\left(x-1 \right)}^2} \right)}, & \lambda^{- 3} <x <1, \cr 1, & x \geqslant 1. \end{cases} $$

\begin{figure} [htp]
\centerline {\includegraphics [width = 6 true cm] {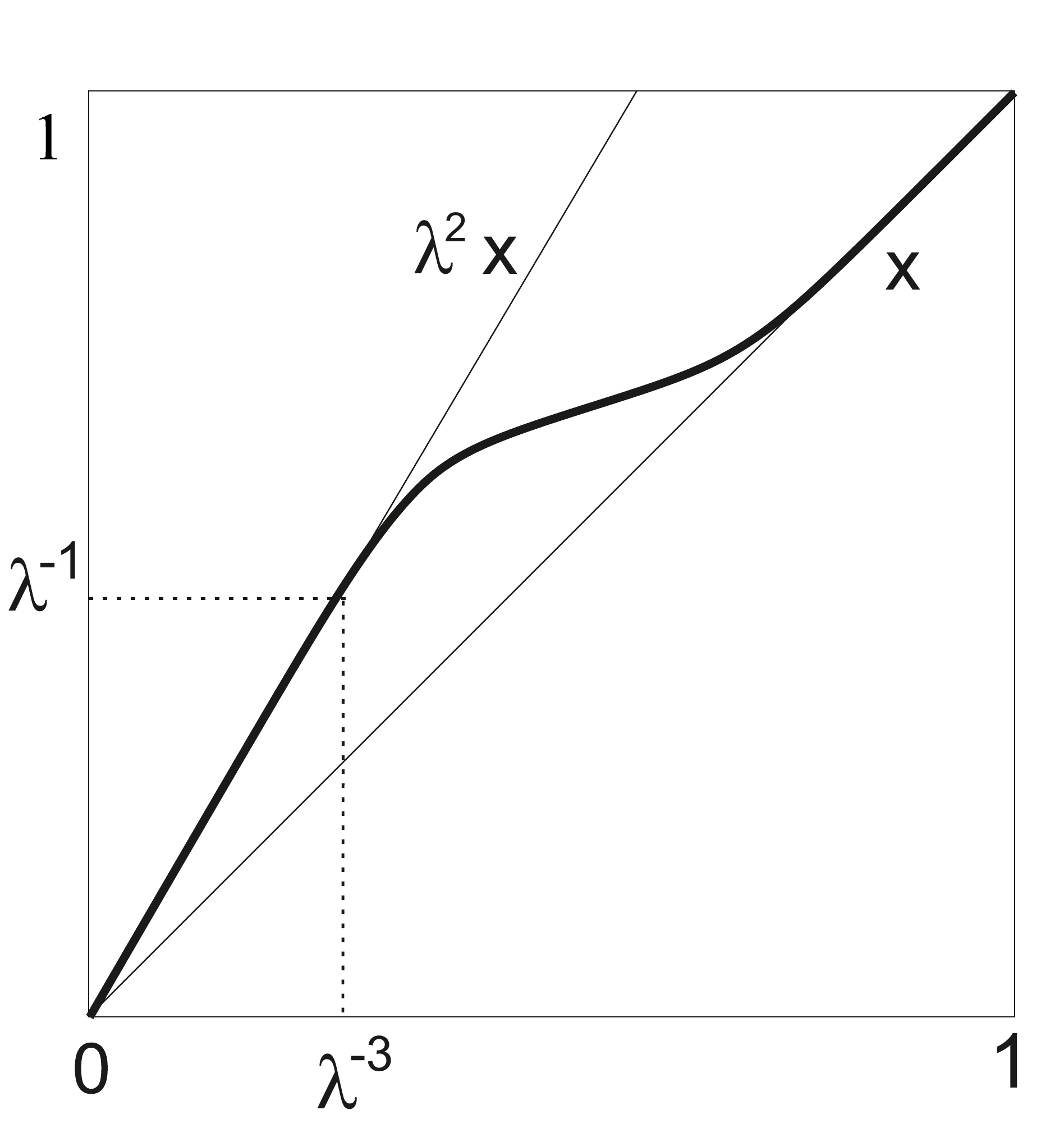}}
\caption{Graph of the function $ \nu (x)$} \label{nu}
\end{figure}

Define a function $ \nu: [0,2] \to [0,2] $ by the formula (see the figure \ref{nu}) $$ \nu (x) = \begin{cases} {\lambda}^2 x , & 0 \leqslant x \leqslant \lambda^{- 3}, \cr \sigma (x) x + (1- \sigma (x)) {\lambda}^2 x, & \lambda^{- 3} < x \leqslant 1,\cr x,&1<x\leqslant 2. \end{cases} $$

Extend the function $ \nu $ on the segment $ [- 2,2] $ by the formula $ \nu (-x) = - \nu (x) $ for $ x \in [0,2] $.

Let $ \mathbb D^2 = \{(x, y) \in \mathbb R^2 | x^2 + y^2 \leqslant 4 \} $. Define a function $ \gamma_A: \mathbb D^2 \to [0,2] $ by the formula
 $$ \gamma_A (x, y) = \begin{cases} \nu (x), & 0 \leqslant | y | \leqslant \lambda^{- 3}, \cr \sigma (| y |) x + (1 - \sigma (| y |) \nu (x), & \lambda^{- 3} <|y| \leqslant 1,\cr x,&1<|y|\leqslant 2. \end{cases} $$
Define a diffeomorphism $ B_A: \mathbb D^2 \to \mathbb D^2 $ by the formula $$ B_A (x, y) = (\gamma_A (x, y), y). $$ By construction, $ B_A (x, y) = (\lambda^2x, y) $ if $ x^2 + y^2 \leqslant \lambda^{- 6} $ and is identical on $ \partial \mathbb D^2 $.

Let $ (x, y) $ be local coordinates in a neighborhood of $U(O)$ of $O$ on $ \mathbb T^2 $ such that the diffeomorphism $ \widehat C $ in these coordinates has the form $$ \widehat C ( x, y) = \left(x / \lambda, \lambda y \right). $$ Then $ Ox \subset W^s_{O} $ and $ Oy \subset W^u_{O} $, as well as $ \{y = const \} $ and $ \{x = const \} $ are stable and unstable foliations. Define a diffeomorphism $ \widehat {B} _A: \mathbb T^2 \to \mathbb T^2 $, which is $ B_A $ in $ U (O) $ and is identical outside $ U (O) $. Then, according to \cite {Wi70}, \cite {KaHa99}, the diffeomorphism $ \widehat \Psi_A = \widehat {B} _A \circ \widehat {C} $ is a $ DA $-diffeomorphism whose non-wandering set consists of a one-dimensional attractor with a unique bunch of degree $2$ and a source fixed point at $O$.

\subsection{Diffeomorphism on the manifold $ \mathbb{T}^2 \times \mathbb{R} $ with a one-dimensional canonically embedded surface attractor} \label{attr}

Consider a smooth function $ \varphi: \mathbb{R} \to \mathbb{R} $ defined by the formula $ \varphi (z) = \frac {z} {\lambda} $. Define a diffeomorphism $ \Phi_A $ on $ \mathbb{T}^2 \times \mathbb{R} $ by the formula $$ \Phi_A (w, z) = (\widehat \Psi_A (w), \varphi (z)). $$

The diffeomorphism $ \Phi_A $ is an $ A $-diffeomorphism whose non-wandering set consists of one fixed saddle point $O\times\{0\}$ and a one-dimensional connected attractor $A$ located on the torus $ \mathbb T^2 \times \{0 \} $. Let us show that the surface attractor $ A $ is canonically embedded.

In $U(O)$ let $ D_1 = \left\{(x, y) ~ | ~ x^2 + y^2 \leqslant \lambda^{- 8} \right\} $ and $ D_2 = \left\{(x, y) ~ | ~ x^2 + y^2 \leqslant \lambda^{- 6} \right\} $. Let $S_A=\mathbb{T}^2 \setminus int \, D_1$ and 
$ U_{A} = S_A\times \left[- \lambda^{- 3}, \lambda^{- 3} \right] $. Then $\Phi_A (U_{A}) = (\mathbb{T}^2 \setminus int \, D_2) \times \left[- \lambda^{- 4}, \lambda^{4} \right] $  and therefore $ \Phi_A (U_{A }) \subset int \, U_{A} $, and $ \bigcap \limits_{n \in \mathbb N} {\Phi_A^n (U_{A})} = A $. This means that $ U_{A}$ is a trapping neighborhood for $ A $. 

\subsection{Diffeomorphism on the manifold $ \mathbb{T}^2 \times \mathbb{R} $ with a one-dimensional canonically embedded surface repeller} \label{repel}
We define a function $ \gamma_R: \mathbb D^2 \to [0,2] $ by the formula $$ \gamma_R (x, y) = \begin{cases} \nu^{- 1} (y), & 0 \leqslant | x | \leqslant \lambda^{- 3}, \cr \sigma (| x |) y + (1- \sigma (| x |) \nu^{- 1} (y), & \lambda^{- 3} <| x | \leqslant 1,\cr y,&1<|x|\leqslant 2. \end{cases} $$
Define a diffeomorphism $ B_R: \mathbb D^2 \to \mathbb D^2 $ by the formula $$ B_R (x, y) = (x, \gamma_R (x, y)). $$ By construction $ B_R (x, y) = \left(x, \frac {y} {\lambda^2} \right) $ if $ x^2 + y^2 \leqslant \lambda^{- 6} $ and is identical on $ \partial \mathbb D^2 $.

Let  $ \widehat {B} _R: \mathbb T^2 \to \mathbb T^2 $ be a  diffeomorphism which is $ {B} _R $ in $U(O)$ and is identical outside $U(O)$. Let $ \widehat \Psi_R = \widehat {B} _R \circ \widehat {C} $. Then, according to \cite{Wi70}, \cite{KaHa99}, the diffeomorphism $ \widehat \Psi_R = \widehat {B} _R \circ \widehat {C} $ is a $ DA $-diffeomorphism whose non-wandering set consists of a one-dimensional repeller $R$ with a unique bunch of degree $ 2 $ and a sink fixed point at $O$.

Consider a copy of the manifold $ \mathbb{T}^2 \times \mathbb{R} $ and a diffeomorphisms $ \Phi_R $ given by the formula $$ {\Phi} _R (w, z) = (\widehat \Psi_R (w ), \varphi^{- 1} (z)). $$
$ {\Phi} _R $ is an $A$-diffeomorphism whose non-wandering set consists of a saddle point $O\times\{0\}$ and a one-dimensional surface canonically embedded repeller $R$ located on the 2-torus $ \mathbb T^2 \times \{0 \} $. The trapping neighborhoods $ U_{R} $ of the repeller $ R $ is exactly the same as for  $A $. Let $ K^{R} _1 = \partial {U_{R}} $, $ K^{R} _2 = {\Phi} _R^{- 1} (K^{R} _1) = K^{R } _2 $.

\section{Proof of theorem \ref{theo0}}
We will realize the following scheme of a diffeomorphism $f$  construction. 
\begin{itemize}
\item consider an Anosov diffeomorphism on two-dimensional torus $ \mathbb{T}^2 $;
\item  perform the Smale ``surgery'' to get a structurally stable $ DA $-diffeomorphism $ \widehat \Psi_A $, whose non-wandering set  consists of a fixed source $ O $ and a one-dimensional attractor (see section \ref{surgery});
\item on the manifold $ \mathbb{T}^2 \times \mathbb{R} $ consider a diffeomorphism $ \Phi_A $ defined as the direct product of the diffeomorphism $ \widehat \Psi_A $ on the torus $ \mathbb{T}^2 $ and a contraction to the origin $ O $ on the line $ \mathbb{R} $. The non-wandering set of the diffeomorphism $ \Phi_A $ consists of the one-dimensional attractor $ A $ and the saddle fixed point (see section \ref{attr});
\item choose a trapping neighborhood $ U_A $ of the attractor $ A $ and the fundamental domain $ K^A = U_A \setminus int \, \Phi_A (U_A) $ of the diffeomorphism $ \Phi_A $ restriction to $ (\mathbb{T}^2 \setminus O) \times \mathbb{R} $. Denote by $ K^A_1 $ and $ K^A_2 = \Phi_A (K^A_1) $ the boundary components of $\partial K^A$;
\item on a copy of the manifold $ \mathbb{T}^2 \times \mathbb{R} $ consider a diffeomorphism
$ \Phi_R $ whose non-wandering set consists of a one-dimensional repeller $ R $ and a saddle fixed point. Moreover, the trapping neighborhood $ U_R $ of the repeller $ R $ coincides with the set $ U_A $, the fundamental domain $ K^R = U_R \setminus int \, \Phi_R^{- 1} (U_R) $ of the  diffeomorphism $ \Phi_R $ restriction to the set $ (\mathbb{T}^2 \setminus O) \times \mathbb{R} $ coincides with the set $ K^A $ and its boundary components $ K^R_1 $ and $ K^R_2 = \Phi^{- 1 } _R (K^R_1) $ coincide with the components $ K^A_1 $ and $ K^A_2 $, respectively (see section \ref{repel});
\item denote by $ M^3 $ a three-dimensional orientable closed manifold which is a result of a  gluing of the fundamental domains by a diffeomorphism $ H: K^{R} \to K^{A} $ with the following properties (see section \ref{gl}):

- $ H (K^{R} _2) = K^A_1 $ and $ H (K^{R} _1) = K^A_2 $;

- $ {\Phi} _A \circ H | _{K^R_2} = H \circ \Phi_R | _{K^R_2} $;
\item denote by $ p: U_A \cup_{H} U_R \to M^3 $ the natural projection. Then the desired diffeomorphism $ f: M^3 \to M^3 $ coincides with the diffeomorphism $ p \Phi_R p^{- 1} | _{p (U_R)} $ on $ p (U_R) $ and with the diffeomorphism $ p \Phi_A p^{- 1} | _{p (U_A)} $ on $ p (U_A) $.
\end{itemize}

Below we give details of the gluing.

\subsection{Gluing of the fundamental domains}\label{gl}
In this section we construct a diffeomorphism $ H: K^{R} \to K^{A} $ with the following properties:
\begin{itemize}
\item $ H (K^{R} _2) = K^A_1 $ and $ H (K^{R} _1) = K^A_2 $;
\item $ {\Phi} _A \circ H | _{K^R_2} = H \circ \Phi_R | _{K^R_2} $.
\end{itemize}

Before the construction we explain how to construct the desired diffeomorphism $ f $ using the diffeomorphism $ H $.

Let $ M^3 = U_A \cup_{H} U_R $ be a factor space obtained by a gluing $ U_A $ and $ U_R $ with respect to $ H $. Denote by $ p: U_A \cup_{H} U_R \to M^3 $ the natural projection and by $ f: M^3 \to M^3 $ a diffeomorphism coinciding with the diffeomorphism $ p \Phi_R p^{-1} | _{p (U_R)} $ on $ p (U_R) $ and with the diffeomorphism $ p \Phi_A p^{- 1} | _{p (U_A)} $ on $ p (U_A) $.
 By construction $ f $ is an $ A $-diffeomorphism whose non-wandering set is the union of connected one-dimensional canonically embedded surface  repeller $ R $ and attractor $ A $.

Now describe the construction of the diffeomorphism $ H $.

\begin{figure} [htp]
\centerline{\includegraphics [width = 9 true cm] {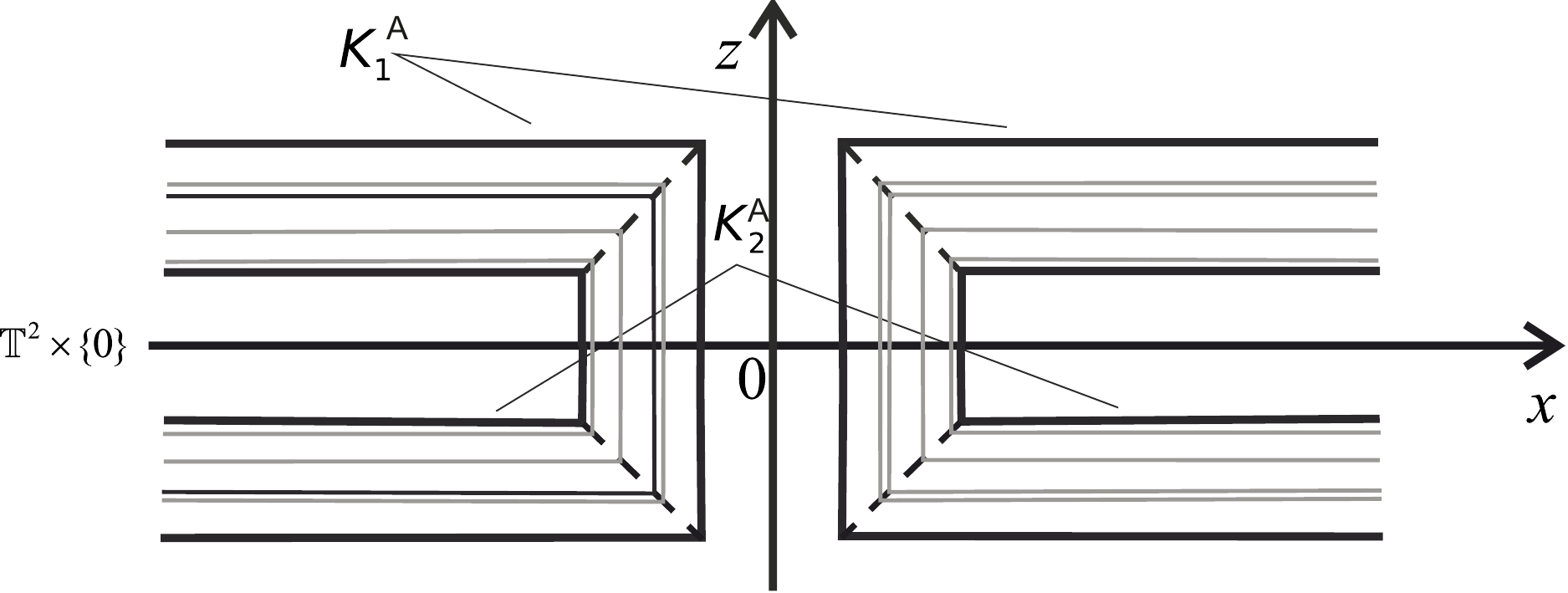}}
\caption{Foliation of the fundamental domain $ K^A $} \label {fdA}
\end{figure}

Let us foliate the fundamental domain $K^A$ by pretzels of genus $2$ in the following way. Define a linear function $ r: [0,1] \to \left[\lambda^{- 4}, \lambda^{- 3} \right] $ by the formula $$ r (t) = (\lambda^{- 3 } - \lambda^{- 4}) t + \lambda^{- 4}. $$
Let $ D_{r (t)} = \{(x, y) \in U (O): x^2 + y^2 \leqslant r^2 (t) \} $ and $ G_t, \, t \in [0,1] $, coincides with the set $ \mathbb T^2 \times \{- r (1-t), r (1- t) \} $ outside of $ D_{r (t)} \times \mathbb R $ and coincides with the cylinder $ \partial D_{r (t)} \times [-r (1-t), r (1-t )] $ otherwise (see the figure \ref{fdA}). Thus, $G_t$ is a pretzels of genus $2$ and $ G_0 = K^A_1 $, $ G_1 = K^A_2 $. 

\begin{figure} [htp]
\centerline {\includegraphics [width = 9 true cm] {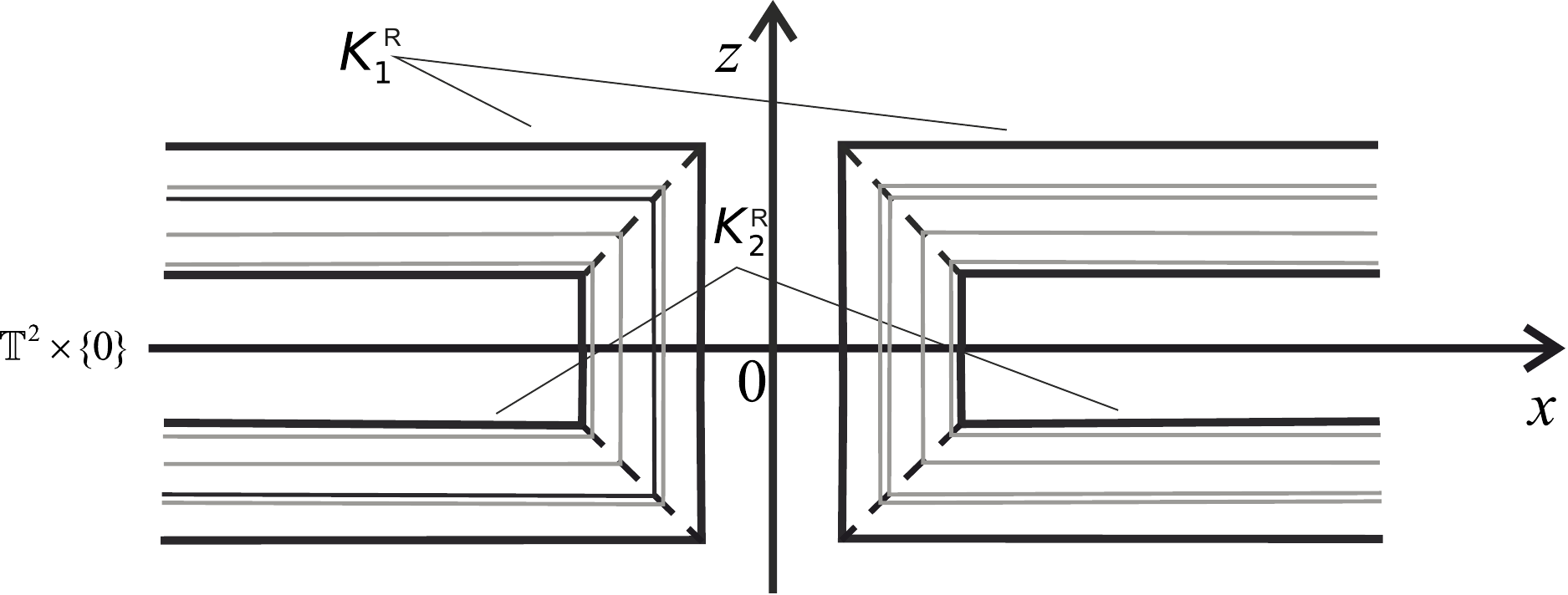}}
\caption {Foliation of the fundamental region $ K^R $} \label {fundamental_domain_foliation}
\end{figure}

As $ K^R $ is a copy of $ K^A $ then $ K^R $ is also foliated by the pretzels $ G_t, \, t \in [0,1] $ so that $ G_0 = K^R_1 $ and $ G_1 = K^R_2 $ (see figure \ref{fundamental_domain_foliation}).

In the local coordinates $ (x, y) $ of the neighborhood $ U (O)$, on the disk $ \mathbb D^2 $ we define the diffeomorphism $ B_t, t \in [0,1] $ by the formula $$ {B} _t (x , y) = t {B} _R (x, y) + (1-t) B_A (x, y). $$ Let $ \widehat {B} _t: \mathbb T^2 \to \mathbb T^2 $  coincides with $ B_t $ in $ U (0) $ and is identical outside $ U (0) $. Let $$ \widehat \Psi_t = \widehat C \circ \widehat B_t. $$
Notice that $ \widehat \Psi_0 = \widehat \Psi_A $ and $ \widehat \Psi_1 = \widehat \Psi_R $. Consider the inverse to $r(t)$  function $ \mu (r): \left[\lambda^{- 4}, \lambda^{- 3} \right] \to [0,1] $ given by the formula $$ \mu (r) = \frac {\lambda^4 r-1} {\lambda-1}. $$
For $t\in [0,1]$ let $ \eta (t) = \frac {t} {\lambda} + (1-t) \lambda. $
We define an orientation-changing diffeomorphism $ \kappa $ of the segment $ \left[\lambda^{- 4}, \lambda^{- 3} \right] $ by the formula $$ \kappa (r) = \eta ((\mu (r) ) r. $$
It is directly verified that the diffeomorphism $ \kappa $ has a unique fixed point $ r _* \in \left[\lambda^{- 4}, \lambda^{- 3} \right] $. For $ t \in [0,1] $ Let $$ \tau (t) = \mu (\kappa (r (t))). $$ By construction, $ \tau $ is an orientation-changing diffeomorphism of the segment $ [ 0,1] $ with a unique fixed point $ t _* \in [0,1] $.

\begin{figure} [h]
\centerline{\includegraphics [width = 9 true cm] {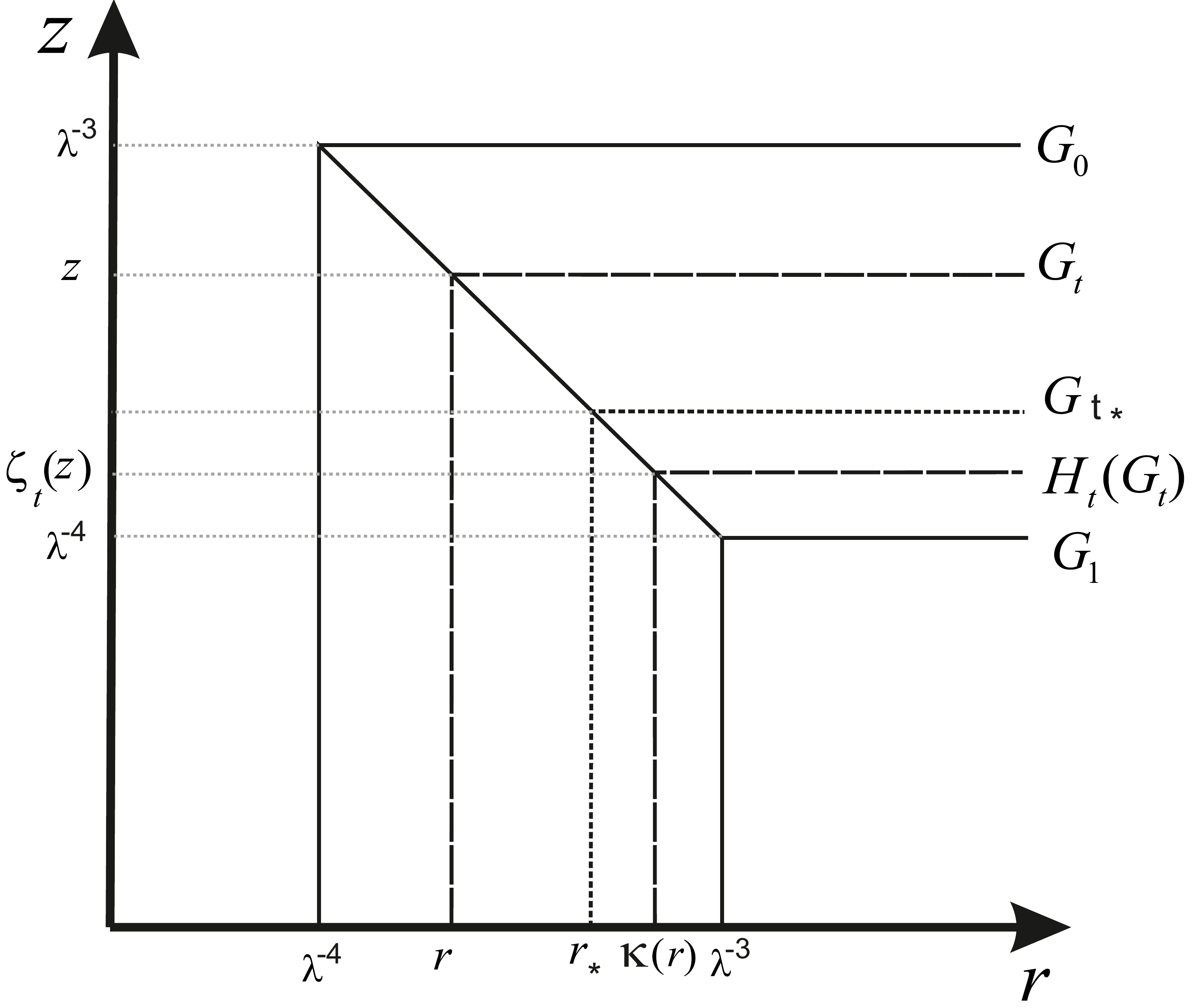}}
\caption{Construction the map $ H_t $} \label{zeta_construction}
\end{figure}

For $ t \in [0,1], \, z \in \mathbb R $ let $ \zeta_t (z) = \frac {r (1- \tau (t))} {r (1-t)} z $. It is directly verified that $ \zeta_{t _*} (z) = z $ for any $ z \in \mathbb R $. On the pretzel $ G_t $ we define $ H_t $ as a diffeomorphism on the image  by the formula $$ H_t (w, z) = (\widehat \Psi_t (w), \zeta_t (z)). $$
It is directly verified that $ H_t (G_t) = G _{\tau (t)} $ (see figure \ref{zeta_construction}).

Denote by $ H: K^R \to K^A $ a diffeomorphism that maps a point $ (z, w) $ of the pretzel $ G_t \subset K^R $ to the point $ H_t (w, z) $ of the pretzel $ G _{\tau (t)} \subset K^A $.

It remains to verify the equality $ {\Phi} _A \circ H | _{K^R_2} = H \circ \Phi_R | _{K^R_2} $, which  is equivalent to the equality $$ {\Phi} _A \circ H_1 (w, z) = H_0 \circ \Phi_R (w, z).$$ 
The last equality follows from  equalities below
$$ \Phi_A \circ H_1 (w, z) = \Phi_A (\widehat \Psi_R (w), \lambda z) = (\widehat \Psi_A \circ \widehat \Psi_R (w), z) $$ and $$ H_0 \circ \Phi_R (w, z) = H_0 (\widehat \Psi_R (w), \lambda z) = (\widehat \Psi_A \circ \widehat \Psi_R (w), z). $$

Note that in this example both the attractor and the repeller lie on a pretzel and resulting diffeomorphism.

\subsection{Behaviour of invariant manifolds of the diffeomorphism $f$}

Denote by $F^s_A$ a foliation on $(\mathbb T^2\setminus O)\times\mathbb R$ consisting of stable manifolds $W^s_a,\, a\in A$. By construction $W^s_a=w^s_a\times\mathbb R$, where $w^s_a\subset\mathbb T^2$ is a stable manifold of the point $a$ with respect to the Anosov diffeomorphism $\widehat C$. The unstable foliation $F^u_R$ looks the same way. 

Thus, the tangent plane to the leaf $W^s_a$ at the point $a=(w,z)\in\mathbb T^2\times\mathbb R$ is stretched over the vectors $ v_1^{s}(w,z), v_2^{s}(w,z)$, where the vector $ v_1^{s}(w,z)$ lies in the plane tangent to the torus $\mathbb T^2\times\{z\}$ and the vector $ v_2^{s}(W, z)$ is perpendicular to it. That is the tangent plane to the leaf $W^s_a$ is parallel to $\mathbb R$, we call it {\it vertical}. Similarly, the tangent plane to the leaf $W^u_r$ at the point $r=(w,z)\in\mathbb T^2\times\mathbb R$ is stretched over the vectors $ v_1^{u}(w, z), v_2^{u}(w, z)$, where the vector $ v_1^{u}(w, z)$ lies in the plane tangent to the torus $\mathbb T^2\times\{z\}$ and the vector $ v_2^{u}(W,z)$ is perpendicular to it. That is the tangent plane to the leaf $W^u_r$ is vertical.

Then the tangent plane to the leaf $H(W^u_r)\cap K^A$ at the point $\bar r=(\bar w,\bar z)=H(r)$ is stretched over the vectors $\bar v_1^{u}(\bar w,\bar z)=DH_{r}(v_1^{u}(w,z)),\bar v_2^{u}(\bar w,\bar z)=DH_{r}(v_2^{u}(w,z))$.  Since the gluing diffeomorphism $H$ preserves the natural splitting $\mathbb T^2\times\mathbb R$ on tori and straight lines, the vector $\bar v_1^{u} (\bar w,\bar z)$ lies in the plane tangent to the torus $\mathbb T^2\times\{\bar z\}$ and the vector $\bar v_2^{u} (\bar w,\bar z)$ is perpendicular to it.

Thus, the two-dimensional foliation $F^s_A$ and $H(W^u_R)$ can have a tangency. Moreover, due to  \cite{RoWi}, the constructed diffeomorphism $f$ has contacts along four planes. That is, the diffeomorphism $f$ is not structurally stable. Moreover, the tangency is not generic.  In the next section, we will show how to tweak the gluing so that the resulting diffeomorphism has a generic tangency.

\subsection{Gluing of the fundamental domains that implies a generic tangency}

\begin{figure}[h]
\centerline{\includegraphics[width=10 true cm]{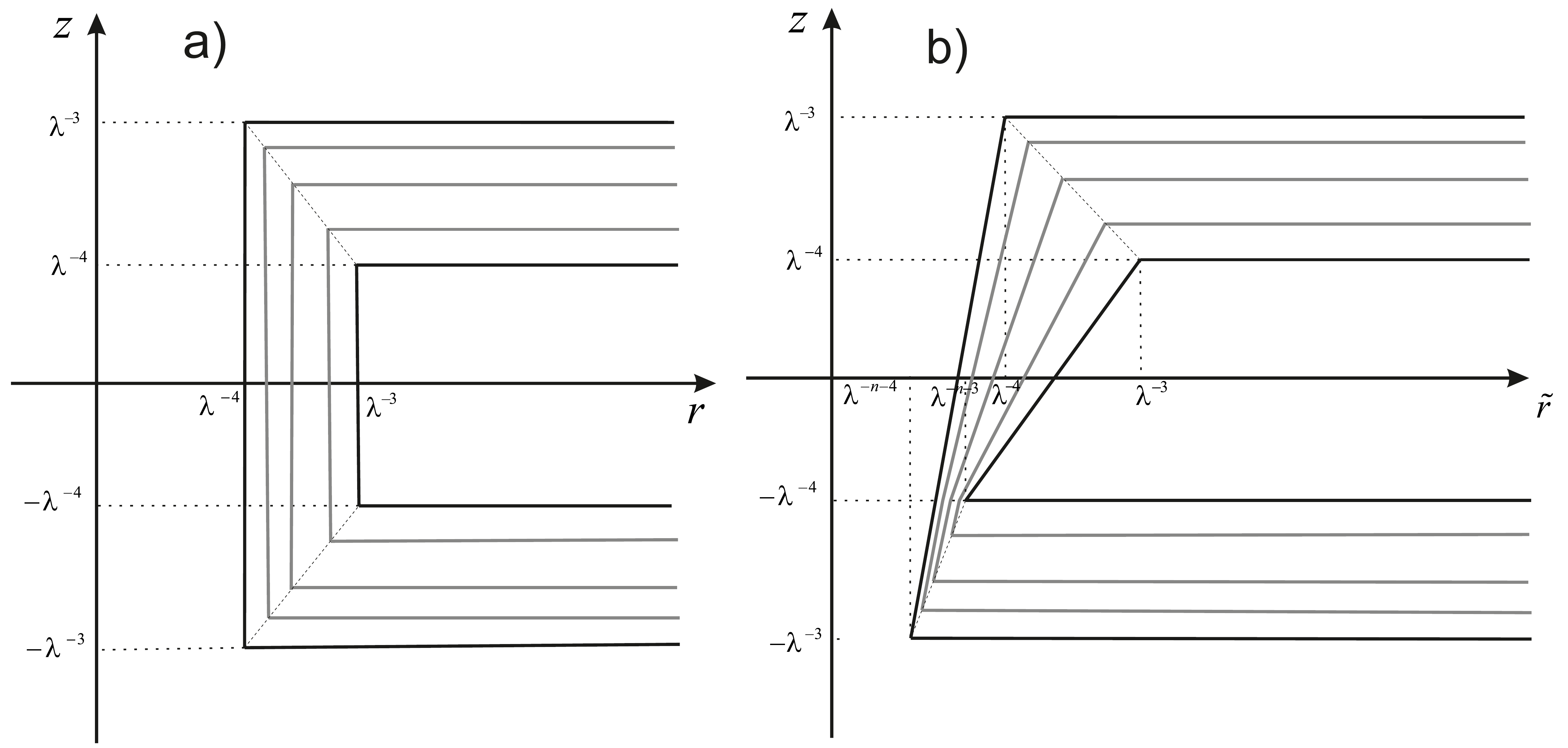}}
\caption{a) Old fundamental domain; b) New fundamental domain}\label{no}
\end{figure}

To construct a new gluing, we will use the previous fundamental domain $K^R$, which is foliated by pretzels $G_t,\, t\in[0,1]$ and the new fundamental domain $\widetilde K^A$ with the new foliation $\widetilde G_t,t\in[0,1]$ (see figure \ref{no}).

We divide each leave $G_t$ in the  fundamental domain $K^A$ of the attractor $A$ into three parts: two tori with a hole $(\mathbb{T}^2\setminus int \,D_{r(t)})\times\{r(1-t)\}$, $(\mathbb{T}^2\setminus int \, D_{r(t)})\times\{-r(1-t)\}$ and the cylinder $\partial D_{r(t)}\times[-r(1-t), r(1-t)]$. Let us construct a map $\theta_t$ which send $G_t$ to a new leaf $\widetilde G_t$. 

Let the coordinates $x$, $y$ in the neighborhood of $U(O)$ are the same as in \ref{surgery}. Firstly, let's construct a map on tori with a hole. Let a map $\theta_t$ be identical on the part $(\mathbb{T}^2\setminus int ~D_{r(t)})\times\{r(1-t)\}$ and $\theta_t=\widehat\Psi^{-n}_A$, for some $n\in\mathbb{N}$ (we can choose $n$ absolutely arbitrary) on $(\mathbb{T}^2\setminus int~D_{r(t)})\times\{-r(1-t)\}$. 
It remains to define $\theta_t$ on cylinders $\partial D_{r(t)}\times[-r(1-t), r(1-t)]$ so that the parts of the leaves are joined. Let's use cylindrical coordinates $r,\varphi,z$. If a  point on the cylinder $\partial D_{r(t)}\times[-r(1-t), r(1-t)]$ has coordinates $r(t), \varphi, z$, then for $t\in\{0,1\}$  its image under the map $\theta_t$ will have coordinates $(\widetilde r(t,z),\varphi, z)$, where $$\widetilde r(t, z) = \frac{\lambda^n-1}{2\lambda^n}\frac{r(t)}{r(1-t)}z+\frac{\lambda^n+1}{2\lambda^n}r(t).$$
If $t\in(0,1)$ then $\widetilde r(t,z)$ is a piecewise linear function given by the following way. Let $\widetilde r(t,-\lambda^{-4})=\left(\widetilde r\left(1,-\lambda^{-4}\right)-\widetilde r\left(0,-\lambda^{-4}\right)\right)t+\widetilde r\left(0,-\lambda^{-4}\right)$ then $$\widetilde r(t,z)=\begin{cases}\frac{\widetilde r\left(t,-\lambda^{-4}\right)-r(t)}{-\lambda^{-4}-r(1-t)}\left(z-r(1-t)\right)+r(t),\cr  z\in[-\lambda^{-4},r(1-t)];\cr\frac{\widetilde r\left(t,-\lambda^{-4}\right)-\frac{1}{\lambda^n}r(t)}{-\lambda^{-4}+r(1-t)}\left(z+r(1-t)\right)+\frac{1}{\lambda^n}r(t),\cr z\in[-r(1-t),-\lambda^{-4}].\end{cases}$$

Denote by $\Theta$ the map composed by $\theta_t$ on $G_t$.  Let's check that $\Theta\circ\Phi_A|_{K^A_1}=\Phi_A\circ\Theta|_{K^A_1}$. Indeed, on an upper torus with a hole, the condition is satisfied because $\Theta$ is identical on it. On a lower torus with a hole, the diffeomorphism $\Theta$ has the form $(\widehat \Psi_A^{-n}, z)$. So for $\forall (w, z)\in (\mathbb{T}^2\setminus int ~D_{r(0)})\times\{-\lambda^{-3}\}$ we have:
$$\Theta\circ\Phi_A(w,z)=\Theta(\widehat\Psi_A(w), z/\lambda)$$ $$=(\widehat\Psi_A^{-n}\circ\widehat\Psi_A(w),z/\lambda)=(\widehat\Psi_A^{1-n}(w),z/\lambda).$$
$$\Phi_A\circ\Theta(w,z)=\Phi_A(\widehat\Psi_A^{-n}(w), z)$$ $$=(\widehat\Psi_A\circ\widehat\Psi_A^{-n}(w),z/\lambda)=(\widehat\Psi_A^{1-n}(w),z/\lambda).$$

On the cylinder $\partial D_{r(0)}\times[-\lambda^{-3},\lambda^{-3}]$ for every point $(r,\varphi,z) \in \partial D_{r(0)}\times[-\lambda^{-3},\lambda^{-3}]$ we have:
$$\Theta\circ\Phi_A(r(0),\varphi,z)=\theta_1(\lambda r(0),\varphi,z/\lambda)=\theta_1(r(1),\varphi,z/\lambda)=$$ $$=\left(\frac{\lambda^n-1}{2\lambda^n}\frac{\lambda^{-3}}{\lambda^{-4}}z/\lambda+\frac{\lambda^n+1}{2\lambda^n}\lambda^{-3},\varphi,z/\lambda\right)$$ $$=\left(\frac{\lambda^n-1}{2\lambda^n} z+\frac{\lambda^n+1}{2\lambda^n}\lambda^{-3},\varphi,z/\lambda\right).$$

$${\Phi}_A\circ \Theta(r(0),\varphi,z)=\Phi_A\circ\theta_0(r(0),\varphi,z)=$$ $$=\Phi_A\left(\frac{\lambda^n-1}{2\lambda^n}\frac{\lambda^{-4}}{\lambda^{-3}}z+\frac{\lambda^n+1}{2\lambda^n}\lambda^{-4},\varphi, z\right)$$ $$=\left(\frac{\lambda^n-1}{2\lambda^{n}}z+\frac{\lambda^n+1}{2\lambda^n}\lambda^{-3},\varphi,z/\lambda\right).$$

Then $\widetilde K^A=\Theta (K^A)$ is a new fundamental domain with a foliation $\widetilde G_t=\Theta(G_t)$, $\widetilde K^A_1 = \Theta (K^A_1)$, $\widetilde K^A_2 = \Theta (K^A_2)$. 

Define the diffeomorphism $\widetilde H:K^{R}\to \widetilde K^{A}$ as $$\widetilde H = \Theta\circ H.$$ Then $\widetilde H$ has properties similar to the diffeomorphism $H$, that is: 

\begin{itemize}
\item $\widetilde H(K^R_2)=\widetilde K^A_1$ and $\widetilde H(K^{R}_1)=\widetilde K^A_2$;
\item ${\Phi}_A\circ \widetilde H|_{K^R_2}=\widetilde H\circ\Phi_R|_{K^R_2}$.
\end{itemize} 

Denote by $\widetilde U_A$ the trapping neighborhood of the attractor $A$ bounded by $\widetilde K^A_1$. Then the ambient manifold as a result of gluing is  the factor space $\widetilde M^3=\widetilde U_A\cup_{\widetilde H}U_R$ and the diffeomorphism $\widetilde f: \widetilde M^3\to \widetilde M^3$  coincids with the diffeomorphism $\widetilde p\Phi_R \widetilde p^{-1}|_{\widetilde p(U_R)}$ on $\widetilde p(U_R)$ and with the diffeomorphism $\widetilde p\Phi_A \widetilde p^{-1}|_{\widetilde p(\widetilde U_A)}$ on $\widetilde p(\widetilde U_A)$, where $\widetilde p:\widetilde U_A\cup U_R\to M^3$ is the natural projection.

Let us prove that the resulting diffeomorphism is structurally stable. To do this, we show that  foliation $F^s_A|_{\tilde K^A}$ is transversal to the foliation $\tilde H (F^u_R)$.

\begin{figure}[h]
\centerline{\includegraphics[width=9 true cm]{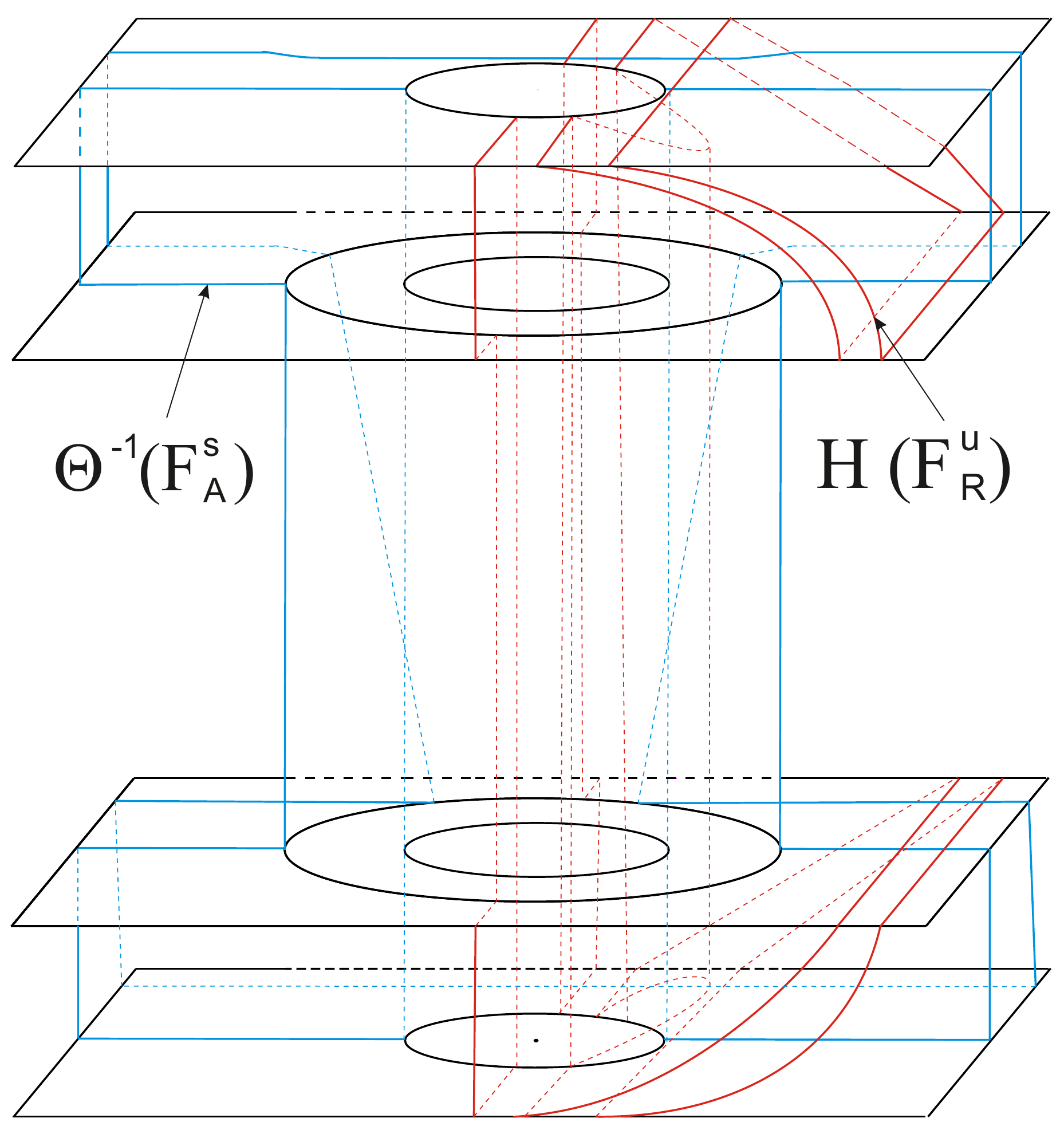}}
\caption{Leaves of the foliation $\Theta^{-1}(F^s_A)$ and $H(F^u_R)$ in $K^A$}\label{fol}
\end{figure}

By construction, transversality is present at all points $\tilde H(w,z)$, where $(w,z)$ belongs to the part of the leaf $G_t$ that is a torus with a hole (see figure \ref{fol}). Therefore, we will check the transversality at the points $\tilde H (w,z)$, where $(w, z)$ belongs to the part of the $G_t$ being the cylinder. At these points, the foliation $F^s_A $  has a form $y=const$ and the foliation $F^u_R$  has a form $x=const$. Note that the tangent planes to leaves of the foliation $H(F^u_R)$ remain vertical, and at points belonging to the plane $y=0$, they have the form $x=const$. 

Since the map $\Theta$ preserves the natural splitting $\mathbb T^2\times\mathbb R$ on tori and is a homothety in $U(O)$ then the tangent planes to leaves of the foliation $\tilde H(F^u_R)$ also is vertical, and at points belonging to the plane $y=0$  have the form $x=const$. Thus, leaves of the foliation $F^s_A$ are  transversal to leaves of the foliation $\tilde H(F^u_R)$ in some neighborhood of the plane $y=0$. At other points, the diffeomorphism $\Theta$ translates any vertical vector at these points to a vector that does not lie in the plane $y=const$ everywhere except the points at the bottom when the foliation $\tilde H(F^u_R)$ is unfolded. The tangency in this case occurs generically along two curves.

\subsection{Why there are infinitely many pairwise $\Omega$-non-conjugated  diffeomorphisms}

Above we constructed examples where the attractor and the repeller belong to a torus with a hole. Obviously, diffeomorphisms constructed from non-similar hyperbolic matrices $C,\,C'$ are not $\Omega$-conjugated.  Moreover, it is clear that the construction admits generalizations to the attractor and repeller obtained by DPA-surgeries that gives infinitely many $\Omega$-non-conjugated  diffeomorphisms.

\section{Proof of Theorem \ref{theo1}} \label{smoothness}

Let $ f$ be an $ A $-diffeomorphism of a closed orientable 3-manifold $M^3$, whose non-wandering set is a union of connected canonically embedded one-dimensional surface attractor $A$ and repeller $ R $. 
By the definition of a canonically embedded attractor, there exists an trapping neighborhood $ U_A $ such that $ U_A \setminus int (f (U_A)) = S_g \times [0,1] = K^A $, and $ K^A $ is the fundamental domain attractor basin. Since the non-wandering set of the diffeomorphism $ f $ consists only of the attractor $ A $ and the repeller $ R $, and they are canonically embedded, the wandering set is representable in the form: $$ M^3 \setminus (A \cup R) = \bigcup \limits_{n = - \infty}^{+ \infty} {f^n (K^A)} = S_g \times (- \infty, + \infty) $$ so that $ \lim \limits_{t \to - \infty} {S_g \times \{t \}} = {R} $ and $ \lim \limits_{t \to + \infty} {S_g \times \{t \}} = {A} $. That is, the wandering set is foliated by the surface $ S_g $.

Then we define the function $ \varphi: M^3 \to [0,1] $ as follows:

$$ \varphi (w) = \left\{\begin{tabular} {ll}
$ \frac {1} {\pi} {\rm arctg} ~ {t} + \frac {1} {2} $, & if $ w \in S_g \times \{t \} $; \\
$ 1 $, & if $ w \in R $; \\
$ 0 $, & if $ w \in A $; \\
\end{tabular} \right. $$

The function $ \varphi $ is a continuous Lyapunov function for $ f $ by the construction. The desired energy function is the result of smoothing the function $ \varphi $ due to the following lemma.

\begin{lemm}\label{l1} Let $ M $ be a smooth compact $ n $-manifold, $ K \subset M $ be a closed subset of $ M $ and $ U $ be some closed neighborhood of $ K $ such that $ K \subset int \, U $. Let a continuous function $ \varphi: {U} \to [0; 1] $ is smooth on $ U \setminus K $ and $ \varphi^{- 1} (0) = K $. Then there exists a $ C^2 $-function $ g: [0; 1] \to [0; 1] $ such that the superposition $ \psi = g \circ \varphi $ is smooth on the whole set $ U $, and the function $ g $ satisfies the following properties:

\begin{itemize}
\item $ g $ monotonically increases on $ [0; 1] $;
\item $ g '(0) = 0 $ and $ g' (c) \neq 0, ~ \forall c \in (0; 1] $;
\item $ g (c) = c, ~ \forall c \in [\frac12; 1] $.
\end{itemize}
\end{lemm}

\begin{demo} Let $ d $ be a metric on the manifold $ M $. For $ c \in (0,1] $ let $ \alpha (c) = \min \{1, d^2 (\varphi^{- 1} (c), K) \} $ and $ \beta ( c) = \max \{1, \max \limits_{x \in \varphi^{- 1} ([c, 1])} ~ | grad ~ \varphi (x) | \} $. By construction the functions $ \alpha (c) $ and $ \beta (c) $ are continuous, and $ \alpha (c) $ is non-decreasing on $ (0; 1] $. Moreover there is a value $ c^* \in (0; 1 ] $ such that $ \alpha (c) $ monotonically increases by $ (0; c^*] $, and $ \beta (c) $ is non-increasing. Then the function $ \gamma (c) = \frac {\alpha (c)} {\beta (c)} $ is a continuous non-decreasing function on the half-interval $ (0; 1] $ and $ \lim \limits_{c \to 0} \frac {\alpha (c) } {\beta (c)} = 0 $.

We construct a $ C^2 $-function $ g: [0; 1] \to [0; 1] $ such that

a) $ g '(c)> 0 $ for any $ c \in (0; 1) $;

b) $ g (c) \leqslant \gamma (c) $ for any $ c \in (0; \frac18) $;

c) $ g '(c) \leqslant \gamma (c) $ for any $ c \in (0; \frac18) $;

d) $ g (c) = c $ for any $ c \in [\frac12; 1] $.

To construct such a function, we will use a partition of unity. Recall that for a given open cover of a topological space $ M $ by open sets $ U _{\alpha} $ with indices $ \alpha $ from the set $ \mathcal A $ {\it is a partition of unity subordinate to the cover $ \{U_\alpha, \alpha \in \mathcal A \} $} is the set of smooth functions $ \sigma_j: M \to \mathbb{R} $, where $ j $ belongs to some set of indices $ J $ with the following properties:

- for every $ j \in J $ there exists $ \alpha \in \mathcal A $ such that $supp (\sigma_j) \subset U_\alpha $, where $supp (\sigma_j) $ is the closure of the set on which function $ \sigma_j $ is non-zero;

- $ 0 \leqslant \sigma_j (x) \leqslant 1 $ for any $ x \in M, j \in J $;

- $ \sum \limits_{j \in J} \sigma_j (x) = 1 $ for any $x \in M$.

If for any point $ x \in M$ there exists a neighborhood $ W_x $ such that the intersection $ W_x \cap supp (\sigma_j) $ is not empty for at most finite number of indices $ j $, then such a partition of unity is called {\it locally finite}.

We take an open cover of the half-interval $ (0; 1] $ by the sets
 
$ U_1 = \{x \in \mathbb{R}: ~ \frac {1} {2} <x \leqslant 1 \} $

$ U_2 = \{x \in \mathbb{R}: ~ \frac {1} {4} <x \leqslant 1 \} $

$ U_i = \{x \in \mathbb{R}: ~ \frac {1} {2^{i}} <x <\frac {1} {2^{i-2}} \} $, $ i = 3.4, \dots $ and the following locally finite partition of unity subordinate to this covering:
$$\sigma_1(x)=\left\{\begin{tabular}{ll}
$1-\sigma_2(x)$, if $x\in(\frac{1}{2};1]$; \\
$0$, if $x\notin(\frac{1}{2};1]$; \\
\end{tabular}\right.$$
$ \forall i = 2,4, \dots $ $$ \sigma_i (x) = \left\{\begin{tabular} {ll}
$ e^{\frac {\left(x- \frac {1} {2^{i-1}} \right)^4} {\left(x- \frac {1} {2^{i}} \right) \left(x- \frac {1} {2^{i-2}} \right)}} $ if $ x \in \left(\frac {1} {2^{i}}, \frac {1} {2^{i-2}} \right) $; \\
$ 0 $ if $ x \notin \left(\frac {1} {2^{i}}, \frac {1} {2^{i-2}} \right) $; \\
\end{tabular} \right. $$
$\forall i=3,5,\dots$ $$\sigma_i(x) =\left\{\begin{tabular}{ll}
$1-\sigma_{i-1}(x)$, if $x\in\left[\frac{1}{2^{i-1}},\frac{1}{2^{i-2}}\right)$; \\
$1-\sigma_{i+1}(x)$, if $x\in\left(\frac{1}{2^{i}},\frac{1}{2^{i-1}}\right)$; \\
$0$, if $x\notin\left(\frac{1}{2^{i}},\frac{1}{2^{i-2}}\right)$; \\
\end{tabular}\right.$$ 
Let $ \varepsilon_i = \gamma \left(\frac {1} {2^{i}} \right) $ for all $ i = 4,5, \dots $. Since each point $ x \in (0,1] $ belongs to the supports of no more than three maps from the partition, the sum $ \sigma_0 (x) = \sum \limits_{i = 4}^{\infty} \varepsilon_i \sigma_i (x) $ is a smooth function on $ (0, \frac12] $. It can be extended continuously by $\sigma_0(0)=0$, similarly to  all functions $ \sigma_i, i \in \mathbb N $. Thus, $ S (x) = \sum \limits_{i = 1}^{\infty} \varepsilon_i \sigma_i (x) $ is a smooth function. Let $ \varepsilon_1 = \varepsilon_2 = 1 $ and $ \varepsilon_3 = \frac {\frac12- \int \limits_0^{\frac12} \sigma_0 (x) dx- \int \limits_0^{\frac12} \sigma_2 (x) dx} {\int \limits_0^{\frac12} {\sigma_3 (x) dx}} $. Define a function $ g $ by the formula $$ g (c) = \left\{\begin{tabular} {ll}
$ \int \limits_{0}^{c} S (x) dx $ if $ c \in (0; 1] $; \\
$ 0 $ if $ c = 0 $ \\
\end{tabular} \right.. $$ It is a $ C^2 $-function, since its derivative is the sum of smooth functions. Let us show that it is the desired one by checking the  conditions a)-d).

a) Since $ g '(c) = S (c) = \sum \limits_{i = 1}^{\infty} \varepsilon_i \sigma_i (c) $, then $ g' (c)> 0 $ for any $ c \in (0; 1) $.

b) For $ i = 4,5, \ldots $ the sequence $ \{\varepsilon_i \} $ is decreasing. Note that for any $ c \in (0; 1] $ there is a unique number $ i^* $ such that $ c \in \left(\frac {1} {2^{i^* - 1}}; \frac {1} {2^{i^* - 2}} \right] $. Then $ \sigma_{i^*} (c) \neq 0 $ and $ \sigma_i (c) = 0 $ for all $ i \notin \{i^*, i^* + 1 \} $. From the choice of the parameters $ \varepsilon_i $ for $ c \in (0, \frac18) $ we obtain the chain of equalities $ g (c) = S (c) = \int \limits_{0}^{c} \left(\sum \limits_{i = 1}^{\infty} \varepsilon_i \sigma_i (x) \right) dx = \int \limits_{0}^{ c} \left(\sum \limits_{i = i^*}^{\infty} \varepsilon_i \sigma_i (x) \right) dx <\int \limits_{0}^{c} \left(\sum \limits_{i = i^*}^{\infty} \varepsilon_{i^*} \sigma_i (x) \right) dx = \varepsilon_{i^*} \int \limits_{0}^{c} \left(\sum \limits_{i = i^*}^{\infty} \sigma_i (x) \right) dx <\varepsilon_{i^*} \int \limits_{0}^{c} \left(\sum \limits_{i = 1}^{\infty} \sigma_i (x) \right) dx = \varepsilon_{i^*} \int \limits_{0}^{c} 1dx = \varepsilon_{i^*} c <\varepsilon_{i^*} = \gamma \left(\frac {1} {2^{i^*}} \right) <\gamma (c) $.

c) For $ g '(c), ~ c \in (0, \frac14) $ the following estimate holds: $ g' (c) = \sum \limits_{i = i^*}^{\infty} \varepsilon_i \sigma_i (c) <\varepsilon_{i^*} \sum \limits_{i = i^*}^{\infty} \sigma_i (c) = \varepsilon_{i^*} <\gamma (c) $.

d) For $ c \in [\frac12; 1] $, the following chain of equalities holds: $ g (c) = \int \limits_{0}^{c} \left(\sum \limits_{i = 1}^{\infty } \varepsilon_i \sigma_i (x) \right) dx = \int \limits_{0}^{\frac12} \left(\sum \limits_{i = 4}^{\infty} \varepsilon_i \sigma_i (x) \right) dx + \int \limits_{0}^{\frac12} \varepsilon_3 \sigma_3 (x) dx + \int \limits_{0}^{\frac12} \varepsilon_2 \sigma_2 (x) dx + \int \limits_{\frac12}^{c} (\varepsilon_1 \sigma_1 (x) + \varepsilon_2 \sigma_2 (x)) dx = \int \limits_{0}^{\frac12} \left(\sum \limits_{i = 4 }^{\infty} \varepsilon_i \sigma_i (x) \right) dx + \varepsilon_3 \int \limits_{0}^{\frac12} \sigma_3 (x) dx + \varepsilon_2 \int \limits_{0}^{ \frac12} \sigma_2 (x) dx + \varepsilon_2 \int \limits _{\frac12}^{c} (\sigma_1 (x) + \sigma_2 (x)) dx = \int \limits_{0}^{\frac12 } \left(\sum \limits_{i = 4}^{\infty} \varepsilon_i \sigma_i (x) \right) dx + \frac {\frac12- \int \limits_0^{\frac12} \left(\sum \limits_{i = 4}^{\infty} \varepsilon_i \sigma_i (x) \right) dx- \int \limits_0^{\frac12} \sigma_2 (x) dx} {\int \limits _0^{\frac12} {\sigma_3 dx}} \int \limits_{0}^{\frac12} \sigma_3 (x) dx + \int \limits_{0}^{\frac12} \sigma_2 (x) dx + (c- \frac12) = c $. So $ g (c) = c $ for $ c \in [\frac12; 1] $.

We show that the superposition $ \psi = g \circ \varphi $ is a smooth function on $ U $.

Notice that $ grad ~ \psi = g '\cdot grad ~ \varphi $, this is useful to us for further discussion. Since the function $ \psi $ on the set $ {U} \setminus K $ is smooth as a superposition of smooth functions, it remains to show that the function $ \psi $ is smooth on the set $ K $.

Consider any point $ a \in K $ and a local chart $ (U_a, h_a) $, where the neighborhood is chosen in such a way that $ \varphi (w) <\frac18 $ for all $ w \in U_a $. First we show differentiability. If the function $ \psi_a = \psi (h_a^{- 1} (x)) $ is differentiable at $ O $, then the function $ \psi $ is differentiable at $ a $. Moreover, the function $ \psi_a $ is differentiable at the point $ O $ and has partial derivatives equal to zero at this point if and only if $ \lim \limits_{s \to O} \frac {\psi_a (s)} {\rho (s, O)} = 0 $, where $ s (x_1, \dots, x_n) \in \mathbb R^n $ and $ \rho $ the Euclidean metric in $ \mathbb{R}^n $, defined by the formula $ \rho (s^1, s^2) = \sqrt {\sum \limits_{i = 1}^n (x_i^1-x_i^2)^2} $ for $ s^1 (x_1^1, \dots, x_n^1), s_2 (x_1^2, \dots, x_n^2) \in \mathbb R^n $. The equality test $ \lim \limits_{s \to O} \frac {\psi_a (s)} {\rho (s, O)} = 0 $ and completes the proof of differentiability.

We introduce a metric $ d_a $ on $ \mathbb R^n $ by the formula $ d_a (s^1, s^2) = d (h^{- 1} _a (s^1), h^{- 1} _a (s^2)) $ for $ s^1, s^2 \in \mathbb{R}^n $. By \cite{Post} (lecture 15), the metrics $ \rho $ and $ d_a $ are equivalent in some compact neighborhood $ U (O) $ of the point $ O $, that is, there exist constants $ 0 <c_1 \leqslant c_2 $ such that what $$ \forall s^1, s^2 \in U (O): ~ c_1 d_a (s^1s^2) \leqslant \rho (s^1, s^2) \leqslant c_2 d_a (s^1, s^2). $$
For $ s \in U (O) $ we put $ w = h_a^{- 1} (s) $ and $ c = \varphi (h_a^{- 1} (s)) = \varphi (w) $. Then $ \lim \limits_{s \to O} \frac {\psi_a (s)} {\rho (s, O)} = \lim \limits_{s \to O} \frac {\psi (h_a^{ -1} (s))} {c_1 d (h_a^{- 1} (s), a)} = \lim \limits_{w \to a} \frac {\psi (w)} {c_1 d (w , a)} = \lim \limits_{w \to a} \frac {g (\varphi (w))} {c_1 d (w, a)} = \lim \limits_{w \to a} \frac { g (c)} {c_1 d (w, a)} <\lim \limits_{w \to a} \frac {\alpha (c)} {\beta (c) c_1 d (w, a)} \leqslant \lim \limits_{w \to a} \frac {d^2 (w, a)} {c_1 d (w, a)} = \lim \limits_{w \to a} \frac {d (w, a )} {c_1} = 0 $.

Now we show that the partial derivatives $ (\psi_a) '_{x_i}, i \in \{1, \dots, n \} $ are continuous at $ O $, that is, $ \lim \limits_{s \to O } (\psi_a) '_{x_i} (s) = 0 $, which is equivalent to $ \lim \limits_{s \to O} | grad ~ \psi_a (s) | = 0 $. Denote by $ J_{h^{- 1} _a} $ the Jacobian of the map $ h_a^{- 1} $, by $ || J_{h^{- 1} _a} || $ its norm subordinated to the Euclidean norm of the vector in $ \mathbb{R}^n $ and by $ B $ a constant such that $ || J_{h^{- 1} _a} (s) || \leqslant B $ for all points $ s $ in some neighborhood of the point $ O $. Then $ \lim \limits_{s \to O} | grad ~ \psi_a (s) | = \lim \limits_{s \to O} | J_{h^{- 1} _a} (s) \cdot g ' (c) \cdot grad ~ \varphi (w) | \leqslant \lim \limits_{s \to O} || J_{h^{- 1} _a} (s) || \cdot | g '(c) | \cdot | grad ~ \varphi (w) | \leqslant \lim \limits_{s \to O} B \cdot \frac {\alpha (c)} {\beta (c)} \cdot | grad ~ \varphi (w) | \leqslant \lim \limits_{w \to a} B \cdot \frac {d^2 (w, a)} {| grad ~ \varphi (w) |} \cdot | grad ~ \varphi ( w) | \leqslant
\lim \limits_{w \to a} B \cdot d^2 (w, a) = 0 $.

Thus, the function $ \psi $ is smooth on $ U $.
\end{demo}



%
%

%

{\bf Acknowledgments.} This work was financially supported by the Russian Science Foundation (project 17-11-01041), except for the \ref{smoothness} section devoted to smoothing the continuous Lyapunov function, written with financial support from the HSE Laboratory of Dynamic Systems and Applications, the Ministry of Science and Higher Education of the Russian Federation, Contract No. 075-15-2019-1931.

\end{document}